\documentclass[12pt]{article}
\usepackage{amssymb,amsmath}
\usepackage{latexsym}
\usepackage{amsfonts}
\newcommand{\qed}{\ifmmode$\Box$\else{\unskip\nobreak\hfil
\penalty50\hskip1em\null\nobreak\hfil$\Box$
\parfillskip=0pt\finalhyphendemerits=0\endgraf}\fi}

\begin{document}
\title
{Spectral factorization in the non--stationary Wiener algebra}
\author
{ D. Alpay, H. Attia, S. Ben--Porat and D. Volok}
\vskip 1truecm
\date{}
\maketitle
\vskip 1truecm
\baselineskip=14pt
\parindent 0cm
\newtheorem{Pa}{Paper}[section]
\newtheorem{Tm}[Pa]{{\bf Theorem}}
\newtheorem{La}[Pa]{{\bf Lemma}}
\newtheorem{Cy}[Pa]{{\bf Corollary}}
\newtheorem{Rk}[Pa]{{\bf Remark}}
\newtheorem{Pn}[Pa]{{\bf Proposition}}
\newtheorem{Pb}[Pa]{{\bf Problem}}
\newtheorem{Dn}[Pa]{{\bf Definition}}
\newtheorem{Ex}[Pa]{{\bf Example}}
\renewcommand{\theequation}{\thesection.\arabic{equation}}

\begin{abstract}
We define the non--stationary analogue of the Wiener algebra
and prove a spectral factorization theorem 
in this algebra. 
%using a factorization theorem for operator--valued
\end{abstract}
%\tableofcontents
\section{Introduction}
\setcounter{equation}{0} In this paper we prove a spectral
factorization theorem in the non--stationary analogue of the
Wiener algebra.
%We then apply this result in combination with the band method
%(see \cite{ggk2}, \cite{gkw-89} and a brief review in Section \ref{damascus}
%for the latter) to solve
%various extension problems in this setting.
To set the problem into perspective and to present our result we
first briefly review the  case of operator--valued functions on
the unit circle. Let ${\mathcal B}$ be a Banach algebra with norm
$\|\cdot\|_{\mathcal B}$ and let ${\mathcal W}({\mathcal B})$
denote the Banach algebra of functions of the form
$$f(e^{it})=\sum_{n\in{\mathbb Z}}e^{int}f_n$$
where the $b_n$ are in ${\mathcal B}$ and such that
$\|f\|_{\mathcal W({\mathcal B})}\stackrel{\rm def.}{=}\sum_{\mathbb Z}
\|f_n\|_{\mathcal B}<\infty$. We set
$$
{\mathcal W}_+({\mathcal B})=\left\{f\in{\mathcal W}({\mathcal
B})\,\,|\,\, f_n=0,\,\,n<0\,\right\}\quad{\rm and}\quad
{\mathcal W}_-({\mathcal B})=\left\{f\in{\mathcal W}({\mathcal
B})\,\,|\,\, f_n=0,\,\, n>0\,\right\}.$$
Inversion theorems in ${\mathcal W}({\mathcal B})$ originate with
the work of Wiener for the case ${\mathcal B}={\mathbb C}$, and
with the work of Bochner and Phillips in the general case; see
\cite{wiener-1932} and \cite{MR4:218g} respectively. Gohberg and Leiterer
studied in \cite{MR48:12113}, \cite{MR48:12114} factorizations in ${\mathcal
W}({\mathcal B})$. A particular case of their results is:
\begin{Tm}
\label{gohlei} Let $W\in{\mathcal W}({\mathcal B})$ and assume
that $W(e^{it})>0$ for every real $t$. Then there exists
$W_+\in{\mathcal W}_+({\mathcal B})$ such that
$W_+^{-1}\in{\mathcal W}_+({\mathcal B})$ and $W=W_+^*W_+$.
\end{Tm}
See Step 3 in the proof of Proposition \ref{8-december-2004}.\\

We recall that $W(e^{it})\in{\mathcal B}$ for every real $t$ and
that $W(e^{it})>0$ is understood in ${\mathcal B}$ (that is,
$W(e^{it})=X(t)^*X(t)$ for some $X(t)\in{\mathcal B}$ which is
invertible in ${\mathcal B}$). For instance, when ${\mathcal B}$
is the space of bounded operators from a Hilbert space into
itself,
$W(e^{it})>0$ means that $W(e^{it})$ is a positive boundedly invertible operator.\\

We note that other settings, where $e^{it}$ is replaced by a
strictly contractive $Z\in \mathcal B$, are also known; see in
particular \cite{MR99k:47035},\cite{ffgk96} and
\cite{MR99i:47027}.\\

In the present paper we consider the case when ${\mathcal B}$ is
the space of block--diagonal operators from $\ell _{{\cal M}}^{2}$
itself (here ${\mathcal M}$ is a pre--assigned Hilbert space and
$\ell_{{\cal M}}^{2}$ denotes the Hilbert space of square summable
sequences with components in ${\cal M}$ and indexed by ${\mathbb
Z}$) and $e^{it}$ is replaced by $Z$, the natural bilateral
backward shift from ${\ell}^2_{{\mathcal M}}$ into itself. This
setting is related to  the theory of non-stationary linear
systems (see \cite{MR99g:93001}). 

\begin{Dn}
The non--stationary Wiener algebra ${\mathcal W}_{NS}$ consists of
the set of operators in ${\mathbf L}({\ell}^2_{{\mathcal M}})$ of
the form $F=\sum_{\mathbb Z} Z^nF_{[n]}$ where the $F_{[n]}$ are
diagonal operators such that
\begin{equation}
\label{mairie-des-lilas} \|F\|_{{\mathcal W}_{NS}}\stackrel{\rm
def.}{=}\sum_{\mathbb Z} \|F_{[n]}\|<\infty
\end{equation}
The element $F$ belongs to ${\mathcal W}_{NS}^+$ (resp.
${\mathcal W}_{NS}^-$) if $F_{[n]}=0$ for $n<0$ (resp. for $n>0$).
\end{Dn}

The {\sl a priori} formal sum $F=\sum_{\mathbb Z} Z^nF_{[n]}$
actually converges  in the operator norm because of
\eqref{mairie-des-lilas}. The fact that ${\mathcal W}_{NS}$ is a
Banach algebra follows from the fact that $ZDZ^*$ and $Z^*DZ$ are
diagonal operators when $D$ is a diagonal operator. \\

The main result of this paper is:

\begin{Tm}
\label{voltaire}
Let $W\in{\mathcal W}_{NS}$ and assume that $W$ positive definite
(as an operator from $\ell^2_{{\mathcal M}}$ into itself).
Then there exists
$W_+\in{\mathcal W}_{NS}^+$ such that $W_+^{-1}\in
{\mathcal W}_{NS}^+$ and $W=W_+^*W_+$.
\end{Tm}

We note the following: as remarked in \cite[p. 369]{MR99g:93001},
Arveson's factorization theorem (see \cite{MR52:3979}) implies
that a positive definite \(W\) can be factorized as $W=U^*U$,
where $U$ and its inverse are upper triangular operators. The new
point here is that if $W$ is in the non--stationary Wiener
algebra, then so are  $U$ and its inverse.\\

A special case of Theorem \ref{voltaire} when $W$ admits a realization
is given in \cite[Theorem 13.5 p. 369]{MR99g:93001} with explicit formula 
for the spectral factor.\\

We now turn to the outline of the paper. It consists of four
sections besides the introduction. In Section \ref{total-recall} we
review the non--stationary (also called time--varying) setting.
In particular we review facts on the Zadeh transform associated
to a bounded upper triangular operator. In Section $3$ we obtain a
spectral factorization for the function $\sum_{{\mathbb Z}}
e^{int}Z^nW_{[n]}$. In Section $4$ we obtain a lower--upper
factorization of that same function. Comparing the two
factorizations lead to the proof of Theorem \ref{voltaire}. This
is done in the last section.

\section{The non--stationary setting and the Zadeh transform}
\label{total-recall} \setcounter{equation}{0} In this section we
review the non--stationary setting. We follow the analysis and
notations of \cite{MR93b:47027} and \cite{DD-ot56}. Let ${\cal
M}$ be a separable Hilbert space, ``the coefficient space''. As
in \cite[Section 1]{DD-ot56}, the set of bounded linear operators
from the space $\ell _{{\cal M}}^{2}$ of square summable sequences
with components in ${\cal M}$ into itself is denoted by ${\cal X}
\left( \ell _{{\cal M}}^{2}\right) $, or ${\cal X}$. The space
$\ell _{{\cal M}}^{2}$ is taken with the standard inner product.
Let $Z$ be the bilateral backward shift operator
\[
\begin{array}{cccc}
\left( Z f\right) _{i}=f_{i+1}, &  & i=\ldots
,-1,0,1,\ldots &
\end{array}
\]
where $f=(\ldots ,f_{-1},
\fbox{$f_{0}$},f_{1},\ldots )\in \ell _{{\cal
M}}^{2}$. The operator $Z$ is unitary on $\ell _{{\cal M}}^{2}$
i.e. $ZZ^{\ast }=Z^\ast Z=I$, and

\[
{\pi}^{\ast}Z^{j}{\pi}= \left\{
\begin{array}{ccc}
I_{{\cal M}} & \mbox{if} & j=0 \\
0_{{\cal M}} & \mbox{if} & j\neq 0.
\end{array}
\right.
\]
where ${\pi}$ denote the injection map
\[
\begin{array}{ccc}
{\pi}:u\in{\cal M}\rightarrow f\in\ell_{{\cal M}}^{2} &
\mbox{where} & \left\{
\begin{array}{cc}
f_{0}=u &  \\
f_{i}=0, & i\neq 0
\end{array}
\right.
\end{array}
.
\]

We define the space of {upper triangular operators} by
\[
{\cal U}\left(\ell_{{\cal M}}^{2}\right) =\left\{A\in%
{\cal X}\left(\ell_{{\cal M}}^{2}\right)
\left|
\begin{array}{ccc}
{\pi}^{\ast}Z^{i}AZ^{\ast j}{\pi}=0 & \mbox{for}
& i>j
\end{array}
\right.\right\},
\]
and the space of {lower triangular operators} by
\[
{\cal L}\left(\ell_{{\cal M}}^{2}\right) =\left\{A\in%
{\cal X}\left(\ell_{{\cal M}}^2\right)
\left|
\begin{array}{ccc}
{\pi}^{\ast}Z^{i} AZ^{\ast j}{\pi}=0 &
\mbox{for} & i<j
\end{array}
\right.\right\}.
\]
The space of diagonal operators ${\cal D}\left(\ell_{{\cal M}}^{2}
\right)$ consists of the operators which are both upper and lower
triangular. As for the space ${\cal X}$, we usually denote these
spaces by ${\cal U}$, ${\cal L}$ and ${\cal D}$.\\

Let $A^{(j)}=Z^{\ast j}AZ^{j}$ for $A\in{\cal X}$ and $%
j=\ldots,-1,0,1,\ldots$; note that
$\left(A^{(j)}\right)_{st}=A_{s-j,t-j}$ and that the maps
$A\mapsto A^{(j)}$ take the spaces ${\cal L}$, ${\cal D}$,
${\cal U}$ into themselves. Clearly, for $A$ and $B$ in
${\mathcal X}$ we have that $\left(AB\right)^{(j)}=A^{(j)}B^{(j)}$ and $%
A^{(j+k)}=\left(A^{(j)}\right)^{(k)}.$\newline

In \cite{MR93b:47027} it is shown that for every $F\in{\cal U}$,
there exists a unique sequence of operators $F_{[j]}\in{\cal D}$,
$j=0,1,\ldots$ such that
\[
F-\sum_{j=0}^{n-1}Z^{j}F_{[j]}\in Z_{{\cal M}}^{n}{\cal U}.
\]
In fact, $\left(F_{[j]}\right)_{ii}=F_{i-j,i}$ and we can
formally represent $F\in{\cal U}$ as the sum of its diagonals
\[
F=\sum_{n=0}^{\infty}Z^{n}F_{[n]}.
\]

More generally one can associate to an element $F\in{\mathcal X}$
a sequence of diagonal operators  such that, formally
$F=\sum_{\mathbb Z}Z^nF_{[n]}$. Recall the well known fact that
even when $F$ is a bounded operator the formal sums
$\sum_{n=0}^{\infty}Z^{n}F_{[n]}$ and
$\sum_{-\infty}^{0}Z^{n}F_{[n]}$ need not define bounded
operators. See e.g. \cite[p. 29]{MR99g:93001} for a counterexample.\\

When the operator $F$ is in the Hilbert--Schmidt class (we will
use the notation $F\in{\mathcal X}_2$) the above representation is
not formal but converges both in operator and Hilbert--Schmidt
norm. Indeed, each of the diagonal operator $F_{[n]}$ is itself a
Hilbert--Schmidt operator and we have:
\begin{equation}
\label{Paul-de-Riquet,1604-1680} \|F\|_{{\mathcal
X}_2}^2=\sum_{n=0}^\infty \|F_{[n]}\|_{{\mathcal X}_2}^2<\infty
\end{equation}
and
\begin{eqnarray*}
\|F-\sum_{-M}^NZ^nF_{[n]}\|^2&\le&
\|F-\sum_{-M}^NZ^nF_{[n]}\|^2_{{\mathcal X}_2} \\
&=&\sum_{-\infty}^{-M-1}\|F_{[n]}\|^2_{{\mathcal X}_2}+
\sum_{N+1}^\infty\|F_{[n]}\|^2_{{\mathcal X}_2}\\
&\rightarrow0&\mbox{{\rm as}}\quad N,M\rightarrow\infty.
\end{eqnarray*}

Here we used the fact that the operator norm is less that the
Hilbert--Schmidt norm:
\begin{equation}
\|F\|\le\|F\|_{{\mathcal X}_2}. \label{Jaures-1859-1914}
\end{equation}
See e.g. \cite[EVT V.52]{MR83k:46003}.

\begin{Dn}
The Hilbert space of upper triangular (resp. diagonal) Hilbert--Schmidt
operators will be denoted by ${\mathcal U}_2$ (resp. by ${\mathcal D}_2$).
\end{Dn}

As already mentioned, the non--stationary Wiener algebra is
another example where the formal power series converges in the
operator norm.

\begin{Pn}
The space ${\mathcal W}_{NS}$ endowed with $\|\cdot\|_{NS}$ is a Banach
algebra.\end{Pn}

{\bf Proof:} Let $F$ and $G$ be in ${\mathcal W}_{NS}$ with
representations
$$F=\sum_{\mathbb Z}Z^nF_{[n]}\quad{\rm and}\quad
G=\sum_{\mathbb Z}Z^nG_{[n]}.$$ Then the family
$Z^mF_{[n]}^{(m-n)}G_{[m-n]}$ is absolutely convergent since
$\|D\|=\|D^{(j)}\|$ for every diagonal operator $D$ and integer
$j\in{\mathbb Z}$.  It is therefore commutatively convergent (see
\cite[Corollaire 1 p. TG IX.37]{tg}) and we can write
\begin{eqnarray*}
FG&=&\left(\sum_{\mathcal Z}Z^nF_{[n]}\right)
\left(\sum_{\mathcal Z}Z^pG_{[p]}\right)\\
&=&\sum_{n,p\in{\mathbb Z}}Z^nF_{[n]}Z^pG_{[p]}\\
&=&\sum_{n,p\in{\mathbb Z}}Z^{n+p}F_{[n]}^{(p)}G_{[p]}\\
&=&\sum_{m\in{\mathbb Z}}Z^m\left(\sum_{n\in{\mathbb
Z}}F_{[n]}^{(m-n)}G_{[m-n]}\right)\\
&=&\sum_{m\in{\mathbb Z}}Z^m(FG)_{[m]}
\end{eqnarray*}
with
$$
(FG)_{[m]}= \sum_{n\in{\mathbb Z}}F_{[n]}^{(m-n)}G_{[m-n]}.$$
This exhibits $FG$ as an element of ${\mathcal W}_{NS}$.
The Banach algebra norm inequality holds in ${\mathcal W}_{NS}$
since
\[
\sum_{m\in{\mathbb Z}}\|(FG)_{[m]}\|\le \sum_{n\in{\mathbb Z}}
\sum_{m\in{\mathbb Z}}
\|F_{[n]}\|\cdot\|G_{[m-n]}\| =\|F\|_{{\mathcal W}_{NS}}\|G\|_{{\mathcal
W}_{NS}}.
\]
\mbox{}\qed

\begin{Dn}
Let $U\in{\mathcal U}$ with formal representation
$U=\sum_{n=0}^\infty Z^nU_{[n]}$. The Zadeh transform of $U$
is the ${\mathcal X}$--valued function defined by
\begin{equation}
U(z)=\sum_{n=0}^\infty z^nZ^nU_{[n]},\qquad z\in{\mathbb D}.
\label{nation}
\end{equation}
\end{Dn}
We note that the series \eqref{nation} converges in the operator norm
for every $z\in{\mathbb D}$ and that $F(z)$ is called in \cite{df1} the symbol
of $F$; see \cite[p. 135]{df1}.
\begin{Tm}
Let $U$, $U_1$ and $U_2$ be upper-triangular operators. Then,
\begin{equation}
\label{place-de-la-republique} \|U(z)\|\le\|U\|,\qquad
z\in{\mathbb D}
\end{equation}
and
\begin{equation}
\label{mardi-25-novembre-2003}
 (U_1U_2)(z)=U_1(z)U_2(z).
\end{equation}
\end{Tm}

{\bf Proof:} A proof of the first claim can be found in
\cite[Theorem 5.5 p 136]{df1}. The key ingredient in the proof is that every
upper triangular contraction is the characteristic function of
a unitary colligation; see \cite[Theorem 5.3 p. 135]{df1}.
To prove the second claim we
remark, as in \cite[p. 136]{df1} that
\begin{equation}
U(z)=\Lambda(z)U\Lambda(z)^{-1}
\label{oberkampf}
\end{equation}
where $z\not =0$ and 
where $\Lambda(z)$ denotes the {\sl unbounded} diagonal operator
defined by
\begin{equation}
\label{Boulets-Montreuil}
\Lambda(z)={\rm diag}~\begin{pmatrix}\cdots& z^{2}I_{\mathcal M}&
zI_{\mathcal M}&I_{\mathcal M}&z^{-1}I_{\mathcal
M}&z^{-2}I_{\mathcal M}&\cdots\end{pmatrix}.
\end{equation}
Of course some care is needed with \eqref{oberkampf}. What is
really meant is that the {\sl a priori} 
unbounded operator on the right coincides
with the bounded operator on the left on a dense set (for
instance on the set of sequences with finite support):
\begin{equation}
\label{samedi-6-decembre-2003}
\Lambda(z)U\Lambda(z)^{-1}u=U(z)u
\end{equation}
where $u\in\ell^2_{{{\mathcal M}}}$ is a sequence with finite
support.\\

We now proceed as follows to prove
\eqref{mardi-25-novembre-2003}. We start with a sequence $u$ as
above. Then:
\begin{align*}
\Lambda(z)^{-1}u&\in\ell^2_{{\mathcal M}}\quad({\rm
since}\,\,u\,\,\mbox{{\rm
has finite support}})\\
U_2\Lambda(z)^{-1}u&\in\ell^2_{{\mathcal M}}\quad({\rm since}\,\,
  {\rm dom}~U_2=\ell_{{\mathcal M}}^2)\\
\Lambda(z) U_2\Lambda(z)^{-1}u&\in\ell^2_{{\mathcal M}}\quad
({\rm by}\,\,\eqref{samedi-6-decembre-2003})\\
\Lambda(z)^{-1}\Lambda(z)U_2\Lambda(z)^{-1}u&\in\ell^2_{{\mathcal
M}}\quad\mbox{{\rm since it is equal to}}~U_2\Lambda(z)^{-1}u\\
\end{align*}
and so (still for sequences with finite support)
$$U_1\Lambda(z)^{-1}\Lambda(z)U_2\Lambda(z)^{-1}u=U_1U_2\Lambda(z)^{-1}u$$
and applying \eqref{samedi-6-decembre-2003} we conclude that
$$
\Lambda(z)U_1\Lambda(z)^{-1}\Lambda(z)U_2\Lambda(z)^{-1}u=
\Lambda(z)F_1F_2\Lambda(z)^{-1}u=(F_1F_2)(z)u\in\ell^2_{{\mathcal
M}}.$$

These same equalities prove \eqref{mardi-25-novembre-2003}.
\mbox{}\qed\mbox{}\\

For a discussion and references on the Zadeh transform we refer to
\cite[p. 255--257]{MR1911847}. The Zadeh--transform was used in 
\cite{MR1878952}, \cite{MR2003f:47021} to attack some problems
where unbounded operators appear. Here our point of view is a bit different.\\

We now extend the Zadeh transform to operators in ${\mathcal
W}_{NS}$ and define
\begin{equation*}
W(e^{it})=\sum_{{\mathbb Z}}e^{int}Z^nW_{[n]}
\end{equation*}
for $W\in{\mathcal W}_{NS}$ with representation $W=\sum_{{\mathbb
Z}}Z^nW_{[n]}$. We obtain a function which is continuous on the
unit circle and it is readily seen that
\[
\|W(e^{it})\|\le\|W\|_{{\mathcal W}_{NS}}\quad
{\rm and}\quad
(W_1W_2)(e^{it})=W_1(e^{it})W_2(e^{it})
\]
for $W,W_1$ and $W_2$ in ${\mathcal W}_{NS}$. To prove the second
equality it suffices to note that $\Lambda(e^{it})$ is now a
unitary operator and that 
\begin{equation}
\label{haiim}
W(e^{it})=
\Lambda(e^{it})W\Lambda(e^{it})^{-1}.
\end{equation}
\section{Spectral factorization of $W(e^{it})$}
\setcounter{equation}{0}
In this section we develop the first step in the proof of Theorem
\ref{voltaire}.
\begin{Pn}
Let $W\in{\mathcal W}_{NS}$ and assume that $W>0$. Then $W(e^{it})>0$ for every
real $t$ and there
exists a ${\mathcal X}$--valued function
$$X(z)=\sum_{n=0}^\infty z^{n}X_n$$
with the following properties:
\begin{enumerate}
\item The $X_n\in{\mathcal X}$ and $\sum_{n=0}^\infty\|X_n\|<\infty$
(that is, $X\in{\mathcal W}_+({\mathcal X})$).

\item $X$ is invertible and its inverse belongs to
${\mathcal W}_+({\mathcal X})$.

\item $W(e^{it})=X(e^{it})^*X(e^{it})$.
\end{enumerate}
\label{8-december-2004}
\end{Pn}
{\bf Proof:} We write $W={\rm Re}~\Phi$ where
$\Phi=W_{[0]}+2\sum_{n=1}^\infty Z^nW_{[n]}$ and proceed in a number of
steps (note that $\Phi$ is a bounded upper triangular operator since
$\sum_{n\ge 0}\|W_{[n]}\|<\infty$).\\

STEP 1. {\sl The operator $(I+\Phi)$ is invertible in ${\mathcal
U}$ and
\begin{equation}
\label{ankara}
(I+\Phi)^{-1}(z)=(I+\Phi(z))^{-1}.
\end{equation}
}

Consider the multiplication operator $M_{\Phi}:{\cal
U}_{2}\longrightarrow{\cal U}_{2}$, defined by $M_{\Phi}F=\Phi F$.
Then,  for every $F,G\in{\cal U}_{2}$, we have
\[\left\langle\Phi F, G\right\rangle_{{\cal U}_{2}}=
\left\langle\Phi F,G\right\rangle_{{\cal X}_{2}}=\left\langle
F,\Phi^* G\right\rangle_{{\cal X}_{2}}. \] Note that \(I+M_\Phi\)
is  a multiplication operator, as well: \(I+M_\Phi=M_\Psi\), where
 $\Psi=\Phi+I$. Hence, for every $F\in{\cal
U}_{2}$, we have:
\begin{eqnarray*}
\left\langle M_{\Psi}F, M_{\Psi}F\right\rangle_{{\cal U}_{2}} & =
& \left\langle\Phi F,\Phi F\right\rangle_{{\cal U}_{2}}+
\left\langle\Phi F,F\right\rangle_{{\cal U}_{2}}+ \left\langle
F,\Phi F\right\rangle_{{\cal U}_{2}}+
\left\langle F,F\right\rangle_{{\cal U}_{2}}\\
& =& \left\langle\Phi F,\Phi F\right\rangle_{{\cal
U}_{2}}+\left\langle\left(\Phi+{\Phi}^{\ast}\right)F,
F\right\rangle_{{\cal X}_{2}}+
\left\langle F,F\right\rangle_{{\cal U}_{2}}\\
& \geq &\left\langle F,F\right\rangle_{{\cal U}_{2}}.
\end{eqnarray*}
In particular, \(M_{\Psi}\) is one-to-one.\\

In the same manner,
\begin{eqnarray*}
\left\langle M_{\Psi}^*F, M_{\Psi}^*F\right\rangle_{{\cal U}_{2}}
& = & \left\langle M_\Phi^* F,M_\Phi^* F\right\rangle_{{\cal
U}_{2}}+ \left\langle M_\Phi^* F,F\right\rangle_{{\cal U}_{2}}+
\left\langle F,M_\Phi^* F\right\rangle_{{\cal U}_{2}}+
\left\langle F,F\right\rangle_{{\cal U}_{2}}\\
& \geq&\left\langle\Phi F,F\right\rangle_{{\cal U}_{2}}+
\left\langle F,\Phi F\right\rangle_{{\cal U}_{2}}+ \left\langle
F,F\right\rangle_{{\cal U}_{2}}\\
& \geq &\left\langle F,F\right\rangle_{{\cal U}_{2}},
\end{eqnarray*}
and, in particular, \(M_{\Psi}\) is onto
(see \cite[p. 30]{Brezis} if need be). Therefore, by the open
mapping theorem $M_\Psi$ is  invertible. \\

Analogous reasoning shows that \(\Psi:\ell^2_{\cal
M}\longrightarrow\ell^2_{\cal M}\) is invertible, as well.
Moreover, the multiplication operator \(M_{\Psi^{-1}}:{\cal
X}_{2}\longrightarrow{\cal X}_{2}\) preserves the subspace \({\cal
U}_{2}\):
\[{M_{\Psi^{-1}}}_{|{\cal U}_{2}}={M_{\Psi^{-1}}}_{|{\cal
U}_{2}}M_\Psi M_{\Psi}^{-1}=M_{\Psi}^{-1}.\] Thus
${\Psi}^{-1}\in{\cal U}$.
Since $(I+\Phi)^{-1}\in{\mathcal U}$ and using
\eqref{mardi-25-novembre-2003} we have:

$$\left((I+\Phi)((I+\Phi)^{-1}\right)(z)=(I+\Phi(z))(I+\Phi)^{-1}(z)$$
and hence we obtain \eqref{ankara}.\qed\mbox{}\\

STEP 2. {\sl It holds that ${\rm Re}~\Phi(z)>0$ for
$z\in{\mathbb D}$.}\\

Indeed, the operator $S=(I+\Phi)^{-1}(I-\Phi)$ is  upper
triangular and $\|S\|<1$. By \eqref{place-de-la-republique},
$\|S(z)\|< 1$ for all $z\in
{\mathbb D}$ and thus ${\rm Re}~(S(z)-I)(S(z)+I)^{-1}>0$.\\

STEP 3. {\sl To conclude it suffices to apply the results of
$\cite{MR47:2412}$ to the function $W(e^{it})$.}\\

Indeed, consider the Toepliz operator with symbol $W(e^{it})$. It
is self--adjoint and invertible since $W(e^{it})>0$. Thus by
\cite[Theorem 0.4 p. 106]{MR47:2412}, 
$W(e^{it})=W_-(e^{it})W_+(e^{it})$ where $W_+$ and its
inverse are in ${\mathcal W}_+({\mathcal X)}$ and $W_-$ and its
inverse are in ${\mathcal W}_-({\mathcal X)}$. By uniqueness of
the factorization $W_+(e^{it})=MW_-(e^{it})^*$. The operator $M$ is
strictly positive and one deduces the factorization result by
replacing $W_+(e^{it})$ by $W_+(e^{it})M^{1/2}$. \qed
\section{Lower-upper factorization of $W(e^{it})$}
\setcounter{equation}{0}
We now use Arveson factorization theorem to obtain another
factorization of $W(e^{it})$.
\begin{Pn}
Let $W\in{\mathcal X}$ be strictly positive. Then there exists
$U\in{\mathcal U}$ such that $W=U^*U$
\end{Pn}
{\bf Proof:} It  suffices to apply the result of \cite{MR52:3979}
(see also \cite[p. 88]{MR84e:93003}) to the nest algebra defined
by the resolution of the identity
$$E_n(\ldots,f_{-1},\fbox{$f_{0}$},f_{1},\ldots )=
(\ldots ,f_{n-1},f_{n})
$$
\mbox{}\qed

\begin{Pn} Let $U,V\in{\mathcal U}$ with formal expansions
$U=\sum_{n=0}^\infty Z^nU_{[n]}$ and $V=\sum_{n=0}^\infty
Z^nV_{[n]}$. Let $z=re^{it}\in{\mathbb D}$. Then
\begin{equation*}
\Omega(z)=V(z)^*U(z)=\sum_{\mathbb Z}e^{imt}Z^m\Omega_{[m]}(r)
\end{equation*}
where
$$
\Omega_{[m]}(r)=
\begin{cases}
&\sum_{p=0}^\infty r^{2p}V_{[p]}^*U_{[p]}\quad {\rm if}\quad m=0 \\
             &\sum_{p=m}^\infty
             r^{2n-m}(V^*_{[p-m]})^{(m)}U_{[p]}\quad{\rm if}\quad m>0
             \end{cases}$$
and $\Omega_{[-m]}=(\Omega_{[m]}^*)^{(-m)}$ and the sums converge
in the operator norm.
\end{Pn}
{\bf Proof:} Since the series $U(z)=\sum_{n=0}^\infty
r^ne^{int}Z^nU_{[n]}$ and $V(z)=\sum_{n=0}^\infty
r^ne^{int}Z^nV_{[n]}$ converge in the operator norm this is an
easy computation which is omitted.

\qed The case $r=1$ is more involved.

\begin{Pn}
\label{steinhaus}
Let $U,V\in{\mathcal U}$ with formal expansions
$U=\sum_{n=0}^\infty Z^nU_{[n]}$ and $V=\sum_{n=0}^\infty
Z^nV_{[n]}$ and let $\Omega=V^*U$. Then the sequence of diagonal
operators associated to $\Omega$ is
$$
\Omega_{[m]}=
\begin{cases}
             &\sum_{p=0}^\infty V_{[p]}^*U_{[p]}\quad {\rm if}\quad m=0\\
             &\sum_{p=m}^\infty (V^*_{[p-m]})^{(m)}U_{[p]}\quad{\rm
             if}m>0.
             \end{cases}$$
where the convergence is entrywise
\end{Pn}

{\bf Proof:} We  first assume that $U$ and $V$ are
Hilbert--Schmidt operator and write
$$U=\sum_{n=0}^NZ^nU_{[n]}+Z^{N+1}R_N\quad{\rm and}\quad
V=\sum_{n=0}^NZ^NV_{[n]}+Z^{N+1}S_N$$ where $R_N,S_N\in{\mathcal
U}$. We denote by $P_0(M)$ the main diagonal of an operator 
$X\in{\mathcal X}$. Then:
$$P_0(\Omega)=\sum_0^{N}V_{[n]}^*U_{[n]}+P_0(S_N^*R_N).$$
The operator $S_N^*R_N$ is Hilbert--Schmidt and so is its main
diagonal $P_0(S_N^*R_N)$. Furthermore by definition of the
Hilbert--Schmidt norm  \eqref{Paul-de-Riquet,1604-1680} and
property \eqref{Jaures-1859-1914} we have
\[
\|P_0(S_N^*R_N)\|_{{\mathcal X}_2}\le\|S_N^*R_N\|_{{\mathcal X}_2}
\le\|S_N^*\|\cdot\|R_N\|_{{\mathcal X}_2}
 \le\|S_N\|_{{\mathcal X}_2}\|R_N\|_{{\mathcal X}_2}
\]
and so
$$\lim_{N\rightarrow
\infty}\|P_0(\Omega)-\sum_{0}^NV_{[n]}^*U_{[n]}\|_{{\mathcal
X}_2}=0.$$
The same holds also in the operator norm thanks to \eqref{Jaures-1859-1914}.\\

Now assume that $U\in{\mathcal U}$. Then we first apply the above
argument to the operators $UD$ and $VE$ where $D,E\in{\mathcal
D}_2$. Then we have for $n=0$
$$E^*\Omega_{[0]}D=\sum_{p=0}^\infty E^*V_{[p]}^*U_{[p]}D$$
where the convergence is in the Hilbert--Schmidt norm. It follows that
$$\Omega_{[0]}=\sum_{p=0}^\infty V_{[p]}^*U_{[p]}$$
where the convergence is entrywise. \mbox{}\qed

\begin{La}
\label{folie-mericourt}
Let $W\in{\cal W}_{NS}$ and let $W=U^*U$ be its Arveson
factorization where $U$ and its inverse are upper triangular. Then
almost everywhere the limit \(U(e^{it}):=\lim_{r\rightarrow 1}
U(re^{it})\) exists in the strong operator topology, has an upper
triangular inverse and satisfies
\begin{equation}
\label{rue-pihet}W(e^{it})=U(e^{it})^*U(e^{it}).
\end{equation}
\end{La}

{\bf Proof:} First, we note that the operator--valued functions
\(U(z),\) \(U^{-1}(z)\),\(U(\overline z)^*\), are analytic in the
open unit disk \(\mathbb D\) and satisfy
\[\|U(z)\|\leq\|U\|,\ \|U^{-1}(z)\|\leq\|U^{-1}\|,
\|U(\overline{z})^*\|\leq\|U\|\] Therefore, by \cite[Theorem A p.
84]{rr-univ} the limits $\lim_{r\rightarrow 1}U^{\pm 1}(re^{it})$,
$\lim_{r\rightarrow 1}U(re^{it})^*$ exist almost everywhere in
the strong operator topology, and  it is easily checked that
\begin{eqnarray*}
\lim_{r\rightarrow 1}U^{-1}(re^{it})&=&U(e^{it})^{-1},\\
\lim_{r\rightarrow 1}U(re^{it})^*&=&\left(\lim_{r\rightarrow 1}
U(re^{it})\right)^{*}=U(e^{it})^*.
\end{eqnarray*}
Note also that if we denote the \(n\)-th diagonal of \(U\) by
\(U_{[n]}\) then for every \(m,n\in\mathbb Z\) and \(F,G\in{\cal
D}_2\) it holds almost everywhere that
\begin{eqnarray*}\langle
U(e^{it})Z^n F,Z^m G\rangle_{{\cal X}_2}& =&\lim_{r\rightarrow
1}\langle U(re^{it})Z^n F,Z^m G\rangle_{{\cal
X}_2}\\
&=&\left\{\begin{array}{c@{\quad }l}\langle e^{i(m-n)t}
Z^{m-n}U_{[m-n]}Z^n F,Z^m G\rangle_{{\cal X}_2},& m\geq n\\
0,& \text{otherwise}\end{array}\right.\end{eqnarray*}
 and hence
\(U(e^{it})\) can be formally written as
\[U(e^{it})=\sum_{n=0}^\infty e^{int}Z^n U_{[n]}.\]

Now, from Proposition \ref{steinhaus} it follows that the 
the diagonals of the operators on both sides of \ref{rue-pihet} coincide
and hence \ref{rue-pihet} follows.
\mbox{}\qed

\begin{Rk} An alternative way to obtain the factorization of $W(e^{it})$
is to define $U(e^{it})=\Lambda(e^{it})U\Lambda(e^{it})^{-1}$. See
equation \eqref{haiim}.
\end{Rk}

\section{Proof of Theorem \ref{voltaire}}
\setcounter{equation}{0}
We proceed in a number of steps:
\begin{enumerate}
\item By Proposition \ref{8-december-2004}
%In the first step we use the result of Section 3:
%we translate the condition $W>0$ into a positivity condition on a
%related function in the Wiener algebra ${\mathcal W}({\mathcal
%X})$ where ${\mathcal X}$ denotes the set of bounded linear
%operators from ${\ell}^2_{{\mathcal M}}$ into itself. We use a
%factorization result of Gohberg and Leiterer to obtain a spectral
we have a factorization $W(e^{it})=X(e^{it})^*X(e^{it})$ where the
${\mathcal X}$--valued function $X(z)$ and its inverse are in
${\mathcal W}_+({\mathcal X})$. At this stage we do not know that
in the series
$$X(z)=\sum_{n=0}^\infty z^nX_n$$
the Fourier coefficients $X_n$ are of the form $X_n=Z^nX_{[n]}$
for some diagonal operator $X_{[n]}$.

\item In the second step we use Arveson's factorization theorem and
Lemma \ref{folie-mericourt}
%start from Arveson's factorization $W=U^*U$
%and consider the lower--upper factorization for $W(e^{it})$
%derived in Section $4$ from Arveson factorisation theorem. More
%precisely
to obtain the factorization
$$W(e^{it})=U(e^{it})^*U(e^{it})$$
where for every real $t$ the operator $U(e^{it})$ and its inverse
are upper triangular. The function $U(e^{it})$ is the limit
function of the Zadeh transform $U(z)$ of $U$; the function
$U(z)$ and its inverse are analytic in ${\mathbb D}$. At this
stage we can write
$$
U(z)=\sum_{n=0}^\infty z^nZ^nU_{[n]}$$ for some diagonal
operators (and similarly for $U^{-1}(z)$), but we do not know
that $U(z)$ (and its inverse) are in ${\mathcal W}_+({\mathcal
X})$, that is whether the sum
$$\sum_{n=0}^\infty\|Z^nU_{[n]}\|=\sum_{n=0}^\infty\|U_{[n]}\|$$
converges or not (and similarly for $U^{-1}(z)$).

\item Comparing the two factorizations derived in Steps 1 and 2 we obtain that
$$X(z)=U(z)M$$
where $M$ is a unitary operator. This leads to
\begin{equation*}
X_nM^*=Z^nU_{[n]},\quad n=0,1,\cdots
\end{equation*}
and allows to obtain the required factorization result for
$W=W(1)$.
\end{enumerate}

\bibliographystyle{plain}
\def\cprime{$'$} \def\cprime{$'$} \def\cprime{$'$}

%\bibliography{/home/user/bib/all}
%\bibliography{/users/faculty/math/dany/bib/all}
%\bibliography{C:/WIN98/Desktop/dany/bib/all}
%\bibliography{/home/alpay/bib/all}
%\bibliography{C:/emacs/dany/bib/all}
%\bibliography{/users/studs/phd/dubi/work/all}
\parbox[t]{8 cm}{Daniel Alpay, Haim Attia, Shirli Ben-Porat and
Dan Volok\\
Department of Mathematics\\
Ben--Gurion University of the Negev\\
Beer-Sheva 84105, Israel\\
{\tt dany@math.bgu.ac.il}\\
{\tt atyah@bgumail.bgu.ac.il}\\
{\tt sbp@math.bgu.ac.il}\\
{\tt volok@math.bgu.ac.il}\\

}
\end{document}